\pdfoutput=1

\documentclass[final,journal]{IEEEtran}

\usepackage[noadjust]{cite}

\usepackage{graphicx}
\graphicspath{{Figures/}}
\DeclareGraphicsExtensions{.pdf,.eps,.png}
\usepackage{changepage}
\usepackage{subfigure}

\usepackage{amsmath}
\usepackage{amssymb}

\usepackage{algorithmic}

\usepackage{float}

\usepackage{url}
\usepackage{hyperref}

\newtheorem{proposition}{Proposition}

\begin{document}
\title{Dynamic Pricing in Shared Mobility on Demand Service}
\author{Han~Qiu,~Ruimin~Li,~and~Jinhua~Zhao% <-this % stops a space
\thanks{}}
% The paper headers
\markboth{IEEE TRANSACTIONS ON INTELLIGENT TRANSPORTATION SYSTEMS}{Qiu \MakeLowercase{\textit{et al.}}: Dynamic Pricing in Shared Mobility on Demand Service}
\maketitle

\begin{abstract}
We consider a profit maximization problem in an urban mobility on-demand service, of which the operator owns a fleet, provides both exclusive and shared trip services, and dynamically determines prices of offers. With knowledge of the traveler preference and the distribution of future trip requests, the operator wants to find the pricing strategy that optimizes the total operating profit of multiple trips during a specific period, namely, a day in this paper. This problem is first formulated and analyzed within the dynamic programming framework, where a general approach combining parametric rollout policy and stochastic optimization is proposed. A discrete-choice-based price optimization model is then used for the request level optimal decision problem and leads to a practical and computationally tractable algorithm for the problem.

Our algorithm is evaluated with a simulated experiment in the urban traffic network in Langfang, China, and it is shown to generate considerably higher profit than naive strategies. Further analysis shows that this method also leads to higher congestion level and lower service capacity in the urban traffic system, which highlights a need for policy interventions that balance the private profit making and the system level optimality.

\end{abstract}


\IEEEpeerreviewmaketitle

\section{Introduction}

\IEEEPARstart{L}{ately}, shared mobility and mobility on-demand (MoD) services have become popular topics in the transportation community. Intuitively, given same amount of vehicle supply, shared mobility has the potential to increase service capacity and reduce congestion. Fagnant \textit{et al.} \cite{Fagnant_2014,Fagnant_2015} used simulation to investigate the extent to which shared autonomous vehicles (SAV), under certain operation rules, can improve transportation efficiencies. Subsequently, they \cite{Fagnant_2015_D} further included dynamic ride sharing in the simulation to show that SAV can also be profitable for the private operator. Alonso \textit{et al.} \cite{Alonso_2017} investigated online matching and dispatching with the view of optimizing the service rate/capacity, and showed that the algorithm could improve the service capacity by a factor of three. However, these studies are overly optimistic as they ignore the private interest of operators by solely focusing on system-level metrics such as service capacity. For example, even though we can optimize operation to increase the capacity with shared mobility, such an operation might be economically sub-optimal for a private operator. This research aims to explicitly model the profit-driven perspective of a private operator and examine the impact of such a profit-driven operation on the performance of the urban traffic system.

Various research efforts have been devoted to optimizing the operations of mobility service operators. There are two major settings: one-sided, or two-sided markets, depending on whether the operator has full control power over the supply. Specifically, one-sided market refers to case in which the operator can manage its fleet in a direct and centralized way and solely transacts with the travelers, while two-sided market refers to the case where the operator is a matching platform and it does business with both drivers and travelers. In this paper, we focus on the one-sided market setting which requires less modeling complexity and admits tractable problem formulations in operation management; for readers interested in the two-sided market setting, we refer them to some recent works in \cite{Banerjee_2016, Chen_2015_u, Jung_2016, Nourinejad_2016}.

Existing studies usually focus on a specific sub-problem instead of formulating the general operation problem integrating dispatching, pricing and rebalancing. Among them, many address the dispatching problem, which is closely related to the traditional vehicle routing problem. One example of the dispatching problem is dial-a-ride problems (DARP). A recent review by Cordeau and Laporte \cite{Cordeau_2007} summarizes attempts on this problem. In terms of revenue management, Atasoy \textit{et al.} \cite{Atasoy_2015} formulated a profit maximization problem with request-level assortment optimization to address the difference between single passenger and shared services. Smith \textit{et al.} \cite{Smith_2013} introduced a fluid model and based on which they developed mathematic programs to optimally rebalance fleet.

This paper presents a formulation of the total profit maximization problem over multiple trips during a day (in contrast to the single request level optimization) of an urban mobility on-demand service. In this problem, we focus on the dynamic pricing strategies of the operator and the corresponding traveler choice behavior. We utilize dynamic programming and choice-based price optimization models \cite{Aydin_2000} to develop a tractable algorithm, and our experiment validates the effectiveness of the proposed algorithm.

The rest of this paper is organized as follows: section II and III present the problem formulation and the algorithm. Section IV evaluates the algorithm in a simulated urban road network in Langfang, China, and examine how the proposed method compares with the myopic pricing strategy and other naive strategies both in terms of the total profit and the impact on traffic congestion and road capacity.

\section{Problem Formulation}

% clearly define the problem

Consider a monopolistic private MoD service operator in an urban transportation market. This operator owns a fleet of (identical) vehicles and provides on-demand travel services. There are two types of service: single or shared. Specifically, each vehicle can provide at most 1 single service at a time, but can serve several shared trips up to its capacity limit, e.g., 3 parties. Vehicle can only provide one type of service at a time, but can change its service type during the course of the day. We consider the problem of finding the optimal pricing strategy of services for each trip request to maximize the total operating profit over a specific period, e.g., a day.
We further make the following assumption to simplify the setting and make the problem tractable:

(1) Traveler or other transportation service provider such as transit agency follow fixed strategies and will not react to different strategies of the operator.

(2) The operator has complete knowledge of the origin-destination distribution of potential trip requests, and the background traffic density (contributed by those not exposed to this MoD service) at each link, at any time. It also has an accurate model of the population traveler choice behavior.

(3) For each vehicle, the operator only controls the assignment of confirmed trips. When providing service, vehicle selects the path with shortest travel time by itself, and then stops at wherever the last passengers disembark and remains idle until the next trip assignment.

Now we formulate our problem formally. Assume during the targeting period (e.g., a day), the operator receives $N$ trip requests $i_1,i_2,\cdots,i_N$ at time $t_1<t_2<\cdots<t_N$. $N$ is a random variable, and each $i_k$ is a random tuple consist of trip origin and destination.
Assume right after the time $t_{k-1}$ of request $i_{k-1}$, the operator has knowledge of current state $Y_k$, which includes current time $t_{k-1}$, fleet status such as the location of each vehicle and details of its assigned trips.
Given $Y_k$, the probability distribution of next request $i_k$ along with the arrival time $t_k$ can be completely determined by the new state $(i_k,t_k)\sim P'(\cdot|Y_k)$.
At time $t_k$ right before receiving the request $i_k$, system state change to $X_k = f'(Y_k,t_k)$: for example, some vehicles move to new location and finish service; some vehicles pick up their passengers. Notation $X$ contains same information as $Y$ is used only to distinguish between the relationships with request $i$: $X$ refers to state before receiving $i$ and $Y$ refers to state after finishing $i$, i.e., confirmation of trip details or rejection of services.

Since not all service options are available to $i_k$ according to $X_k$, we therefore define the feasible set $C_k(i_k,X_k)$ as the collection of feasible service along with the corresponding trip plan. For example, if both single service and shared service are available, $C_k$ will contain two elements corresponding to the trip details of single and shared services respectively. One should note that, even there might be several vehicles available for the same service, in the feasible set we only consider one, which is the best according to some predefined operation rules.

Given $C_k$, the operator tries to generate offer prices $u_k$ for $C_k$ based on a policy $\pi$, current request $i_k$, and current state $X_k$ with $u_k = \pi(i_k,X_k) \in R^{|C_k(i_k,X_k)|}$. Then, the traveler will choose among the feasible services in $C_k$, or reject the offer and consider alternatives, such as driving, taking taxi, or cancelation of the trip. This selection is characterized by a (probabilistic) choice model $w_k\sim P(\cdot|u_k,C_k)$.
For a realization of choice $w_k\in C_k$ or $w_k = reject$, we can calculate the operating profit from this request with $g(u_k,w_k)$, and update the state after request with $Y_{k+1}=f(X_k,w_k)$. Our objective is to maximize the sum of operating profit
\begin{equation}\label{daily_basic}
  \begin{aligned}
  J^* & = \max_{\pi} J(\pi) \\
  & = \max_{\pi} E_{i_1,t_1,w_1,\cdots,i_N,t_N,w_N} \{ \sum_{k=1}^N g(\pi(i_k,X_k),w_k) \}.
  \end{aligned}
\end{equation}

\begin{figure}[!t]
\centerline{\includegraphics[width = 0.5\textwidth]{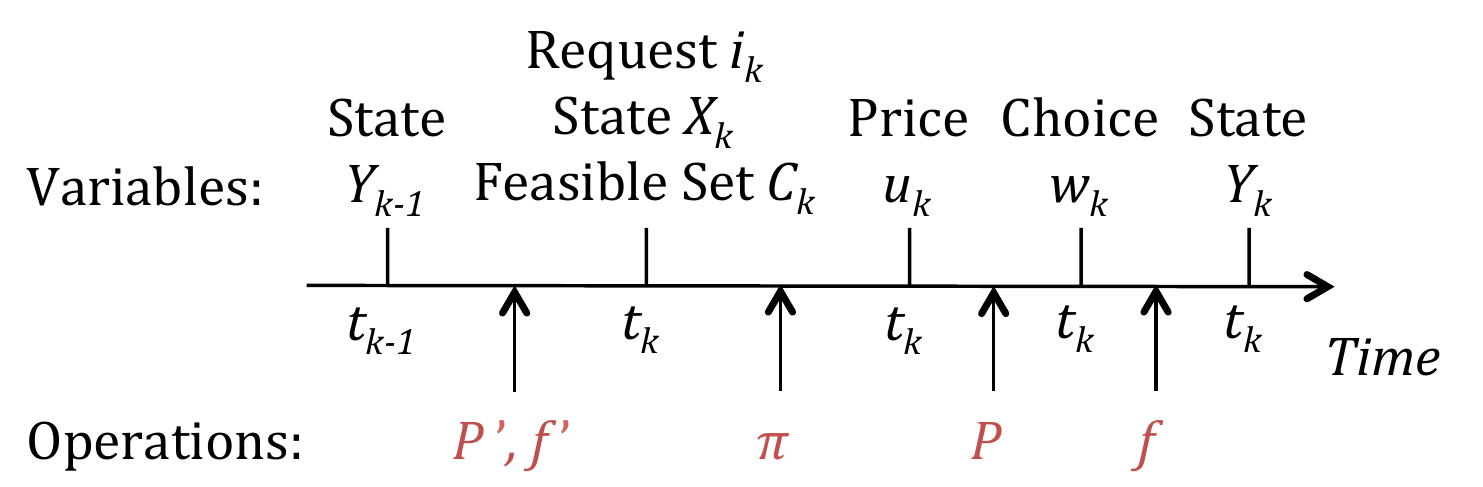}}
\caption{Request Level Dynamics}\label{Fig:flow_dynamic}
\end{figure}

\begin{figure}[!t]
\centerline{\includegraphics[width = 0.4\textwidth]{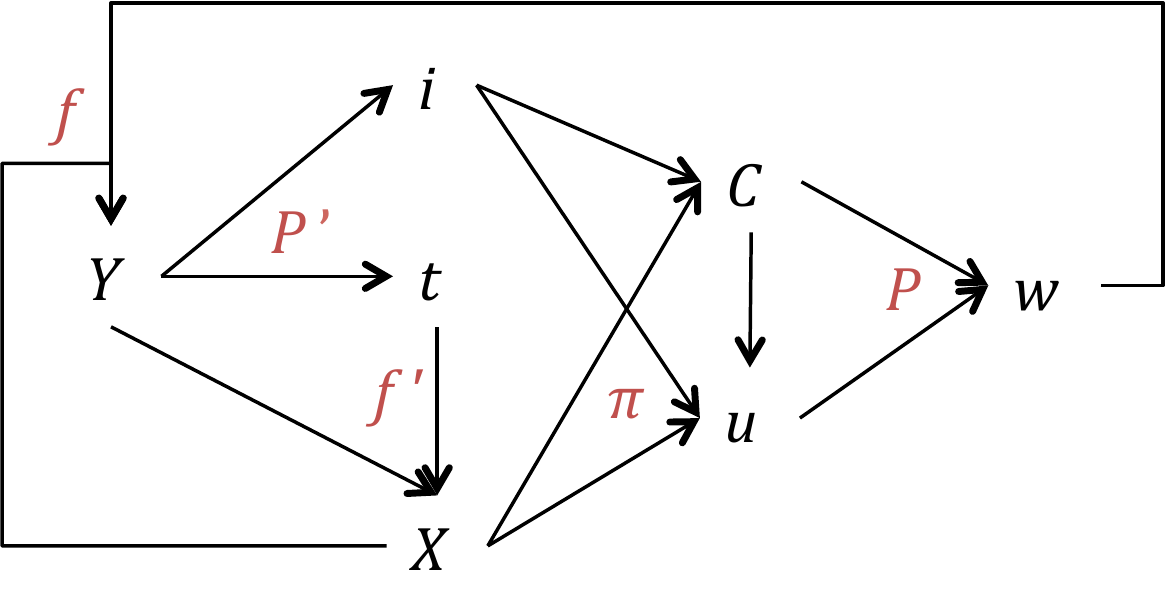}}
\caption{Dependency of Variables}\label{Fig:dependency}
\end{figure}

For reader's convenience, we show the dynamic through time in the Figure \ref{Fig:flow_dynamic}, and the dependency of variables in Figure \ref{Fig:dependency}.

\section{Algorithms}

To arrive at tractable solutions, we first unroll the equation \eqref{daily_basic} into request-wise form, where the process of each request is modeled by the loop $Y\to i,t \to X \to u\to w\to Y$. If we use $J^*(Y)$ to represent the optimal expected total profit in the future given state $Y$, we can write down the following Bellman Equation \cite{Bertsekas_1995}, which states that current expected total profit consists of expected request level profit for next request, and future expected total profit
\begin{equation}\nonumber
  J^*(Y) = E_{i,t|Y}\{\max_{u\in U(C)} E_{w|u,C}\{g(u,w) + J^*(f(X,w))\}\}.
\end{equation}

\noindent where $X=f'(Y,t)$. The optimal, non-myopic request-wise decision for a specific request $i$ at time $t$is then given by
\begin{equation}\label{daily_opt_control}
  u^* = \arg\max_{u\in U(C)} E_{w|u,C}\{g(u,w) + J^*(f(X,w))\}.
\end{equation}

Thus, if the request level subproblem
\begin{equation}\label{request_sub}
  \max_{u\in U(C)} E_{w|u,C}\{g(u,w)\}
\end{equation}

\noindent given $C$ is tractable, and if $J^*(f(X,w))$ can be obtained in a tractable way, the optimal daily-level policy is also tractable. In the following, we discuss the solution of \eqref{daily_basic} when subproblem \eqref{request_sub} can be solve efficiently by some subroutine; we defer the discussion on tractable modeling of \eqref{request_sub} to the next subsection.

Such sequentially decision-making problem is common in transportation research community; to name a few, there are the pricing problem under vehicle routing problem framework \cite{Figliozzi_2007}, the dial-a-ride problems (DARP) \cite{Hyytia_2012}, and the joint dispatching and pricing problem \cite{Sayarshad_2015}.

However, the optimal function $J^*$ is difficult to find in practice, and problem-specific approximation approaches are developed in all studies. Sayarshad and Chow \cite{Sayarshad_2015} provide an exhaustive summary of recent attempts.
In this paper, rather than to approximate optimal value function $J^*$, we conduct approximation directly in the policy space. We notice that, each parametric approximation $\hat{J}(\cdot|\theta)$ with parameter $\theta$, combined with the request level subproblem, provides us a rollout policy $\hat{\pi}_{\theta}$
\begin{equation}\nonumber
  \hat{\pi}_{\theta}(i,X) = \arg\max_{u\in U(C)} E_{w|u,C}\{g(u,w) + \hat{J}(f(X,w)|\theta)\}.
\end{equation}

We can therefore evaluate $\hat{\pi}_{\theta}$ with simulation and optimize $\theta$ with any reasonable optimization methods, such as those of stochastic optimization. This approach is tractable and comparably easier to implement; however, it restricts the policy space and in general we cannot know if this approximation will lead to great gap in performance. Moreover, when the parameter space is large, the optimization become fairly difficult as the evaluation function is indeed a blackbox.

To resolve such potential problems, we introduce an efficient parametric method by utilizing a special structure of our problem. Since $w\sim P(\cdot|u,C)$ relates to the realization of a selection, the request level optimal control problem \eqref{daily_opt_control} can further decomposed into option-level. If, for each trip option $o$ in $C$, we use $\tilde{J}^*(o,X)$ and $p(o|u,C)$ to represent its expected value function and selection probability respectively, and $\tilde{J}^*(0,X)$, $p(0|u,C)$ for rejection, we can write \eqref{daily_opt_control} as follows
\begin{equation}\nonumber
  \begin{aligned}
  u = & \arg\max_{u\in U(C)} \{\sum_{o\in C} p(o|u,C)[g(u,o) + \tilde{J}^*(o,X)] \\
  & + p(0|u,C)\tilde{J}^*(0,X)\} \\
  = & \arg\max_{u\in U(C)} \{\sum_{o\in C} p(o|u,C)[g(u,o) + \tilde{J}^*(o,X)] \\
  & + (1 - \sum_{o\in C} p(o|u,C))\tilde{J}^*(0,X)\} \\
  = & \arg\max_{u\in U(C)} \sum_{o\in C} p(o|u,C)[g(u,o) + \tilde{J}^*(o,X) \\
  & - \tilde{J}^*(0,X)]
  \end{aligned}
\end{equation}

Now, we can further define opportunity cost function $\tilde{A}^*$ to replace the expected value function $\tilde{J}^*$:
\begin{equation}\nonumber
  \tilde{A}^*(o,X) = \tilde{J}^*(0,X) - \tilde{J}^*(o,X), \forall o\in C.
\end{equation}

$\tilde{A}^*$ evaluates the difference of only a single choice selection and it is therefore feasible to approximate it by local information, which is certainly of lower dimension than the global information which is required to describe $\tilde{J}^*$.

We now summarize our approach: we consider the following parametric representation $\hat{A}(o,\hat{X}|\theta)$ with parameter $\theta = \{\theta_s,\theta_{sh}\}$
\begin{equation}\nonumber
  \hat{A}(o,\hat{X}|\theta) =
  \begin{cases}
    \theta_s^T\hat{X} & o \text{ is of type } single\\
    \theta_{sh}^T\hat{X} & o \text{ is of type } shared\\
  \end{cases}
\end{equation}

\noindent along with rollout policy
\begin{equation}\label{control}
  \pi_{\theta}(i,X) = \arg\max_{u\in U(C)} \sum_{o\in C} p(o|u,C)[g(u,o) - \hat{A}(o,\hat{X}|\theta)]
\end{equation}

\noindent where $\hat{X} = (\lambda_{ot},\lambda_{dt},s_{ot},1/v_{ot},1/v_{dt})$ is a local feature tuple including local outflow demand intensity $\lambda_{ot}$ at origin $o$, local outflow demand intensity $\lambda_{dt}$ at destination $d$, local supply capacity of the operator at origin $s_{ot}$, and local average travel speed at both origin and destination $v_{ot},v_{dt}$.
These local values are the average of corresponding variable values at all nodes within some radius $r$ around the target node $o$ or $d$. In this paper, we select $r$ to be $2 \ km$.

As our representation of policy is now low-dimensional, we select an episodic-updated stochastic optimization method to optimize $\theta$: for a given $\theta$, we simulate an episode (i.e., a day) to evaluate the corresponding policy $\pi_{\theta}$ with $J(\pi_{\theta})$, and then use the evaluation to update $\theta$ by Covariance Matrix Adaptation Evolution Strategy (CMA-ES) method \cite{Hansen_2001,Hansen_2003}.

\subsection{Request Level Subproblem}

To complete our algorithm, in this subsection we introduce a tractable model for the request level subproblem. First, let us recall that subproblem \eqref{request_sub}, in option-level, is of the form
\begin{equation}\nonumber
  \max_{u\in U(C)} \sum_{o\in C} p(o|u,C)g(u,o)
\end{equation}

Next, we will introduce price optimization problem under multinomial logit choice model, and show that this problem is suitable and tractable in our setting. For notation simplicity, in the following discussion we assume $C=\{o_1,\cdots,o_n\}$ and use shorthand $j$ to represent option $o_j$; the corresponding price is denoted as $u_j$. The formulation above is therefore simplified as
\begin{equation}\nonumber
  \max_{u\in R^n} \sum_{j=1}^n p(j|u_j)g(u_j,j)
\end{equation}

Assume for each option $j$, there are several values of interest to the traveler and the operator: service type $y_j\in\{single,shared\}$, standard fare $f_j\in R_+$, operating cost $c_j\in R_+$, fare adjustment $\delta_j\in R$, and travel time, i.e., the difference between the estimated arrival time and the current time, $T_j\in R_+$. We first clarify the variable from the constants: rather than considering $u_j = f_j + \delta_j$, we let $u_j = \delta_j$ and only consider $\delta_j$ as variable, while all other are given constant.

Next, we develop the exact form of both $g(u_j,j)$ and $p(j|u_j)$. We let $g(u_j,j)$ to be the profit $g(u_j,j) = R_j = f_j + \delta_j - c_j$, and model $p(j|u_j)$ with the multinomial logit model (MNL) \cite{Luce_1959,Mcfadden_1973,Ben-Akiva_1985}:
\begin{equation}\nonumber
  p(j|u_j) = \frac{e^{U_j}}{\sum_{k=1}^{n} e^{U_k} + e^{U_0}}
\end{equation}

\noindent where the utility $U_j$ of the option $j$ is in form $U_{y_j}(f_j,\delta_j,T_j)$, and $U_0$ is the utility for rejection and given as constant.

With all these details, the request level problem reduces to
\begin{equation}\label{p_p}
  \max_{\delta\in R^n} \frac{\sum_{j=1}^n e^{U_j(\delta_j)}(f_j-c_j+\delta_j)}{\sum_{j=1}^n e^{U_j(\delta_j)} + e^{U_0}}
\end{equation}

This problem is tractable and has good properties when the utility function $U_y$ with respect to variable $\delta$ is continuous and monotonically decreasing for each service type $y$. We first review an important existing proposition:

\begin{proposition}\label{thm:p1} (\cite{Akccay_2010,Gallego_2014,Gallego_2014_b})
  If the utility function $U_y$ with respect to variable $\delta$ is continuous and monotonically decreasing for each service type $y$, to solve problem \eqref{p_p} is equivalent to solve the following equation in $z$
  \begin{equation}\label{e_p}
    e^{U_0}z = \max_{\delta\in R^n} \sum_{j=1}^n e^{U_j(\delta_j)}(f_j - c_j - z + \delta_j)
  \end{equation}
\end{proposition}

The proof can be found in the appendix. Proposition \ref{thm:p1} admits an easy algorithm for a class of utility function specification:

\begin{proposition}\label{thm:p2}
  For the problem \eqref{p_p} with utility $U_j = U_{j0} + E_{y_j}(\delta_j)$ where $U_{j0} = U'_{y_j}(f_j,T_j)$ is constant for all $j$, and $E_{y}(x) = - e_{y1}\min\{0,x\} - e_{y2}\max\{0,x\}$ with $e_{y2} > e_{y1} > 0$ for all service type $y$, the optimal objective value $z^*$ can be solved by the following iterative algorithm:
  \begin{equation}\label{re_p:1}
    \delta_j^{k+1} =
    \begin{cases}
      z^k - f_j + c_j + 1/e_{y_j1} & z^k - f_j + c_j + 1/e_{y_j1} < 0\\
      z^k - f_j + c_j + 1/e_{y_j2} & z^k - f_j + c_j + 1/e_{y_j2} > 0\\
      0 & otherwise
    \end{cases}
  \end{equation}
  \begin{equation}\label{re_p:2}
    z^{k+1} = \frac{\sum_{j=1}^n V_{j0}e^{E_{y_j}(\delta_j^{k+1})}(f_j-c_j+\delta_j^{k+1})}{\sum_{j=1}^n V_{j0}e^{E_{y_j}(\delta_j^{k+1})} + V_0}
  \end{equation}

  \noindent where $V_{j0} = e^{U_{j0}}, V_0 = e^{U_0}$.
\end{proposition}

\section{Simulated Experiment}

We conduct a simulated experiment in a medium-sized city Langfang in China to assess the performance of our algorithm. We first use existing trajectory data of Langfang to calibrate the agent-based simulation system, and then generate different scenarios for validation of the algorithm. In the following subsections, we briefly discuss the implementation of simulation system, and explain the results.

\subsection{Experiment Setup}

\begin{figure}[!t]
\centerline{\includegraphics[width = 0.5\textwidth]{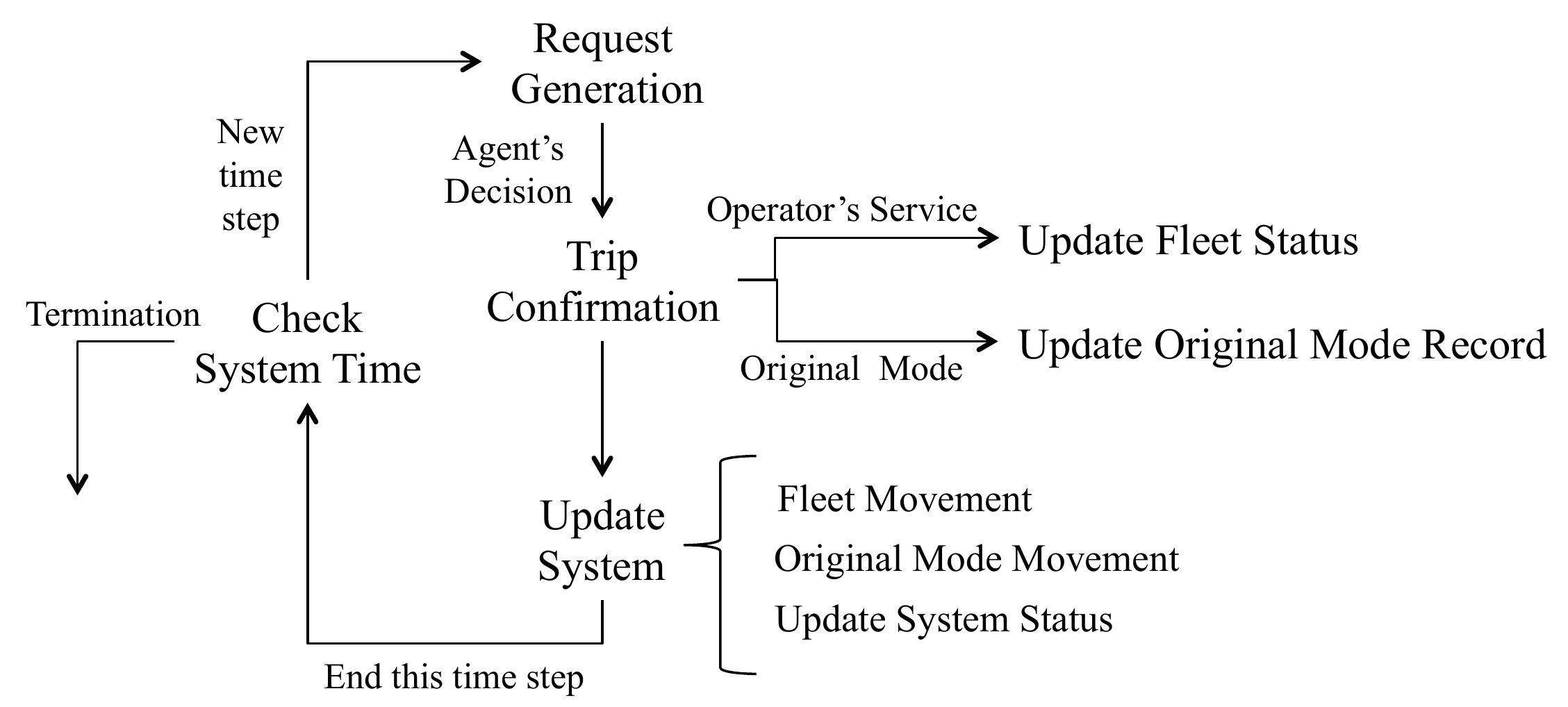}}
\caption{Simulation System Flowchart}\label{Fig:Sim_flowchart}
\end{figure}

To simplify the simulation flow, we discretize the time into minutes and simulate by time steps. At each time step, new demand for the MoD service is generated based on an underlying demand distribution, processed by the operator, and the resulted new passenger trips are stored into the database. After all demand are processed, the system retrieve the alive passenger trips from the database and update the location of corresponding vehicles. The system then ends this time step by updating the network status with an internal speed-density model and the background traffic density. The structure of our simulation system can be described by the flowchart in Figure \ref{Fig:Sim_flowchart}. To ensure efficient training of our algorithm, we further implement following simplification in our simulation system:

(I) We consider driving and taking a taxi as the only travel models and ignore public transportation, walking and cycling to reduce modeling complexity.

\begin{figure}[!t]
  \centering
  \includegraphics[width = 0.5\textwidth]{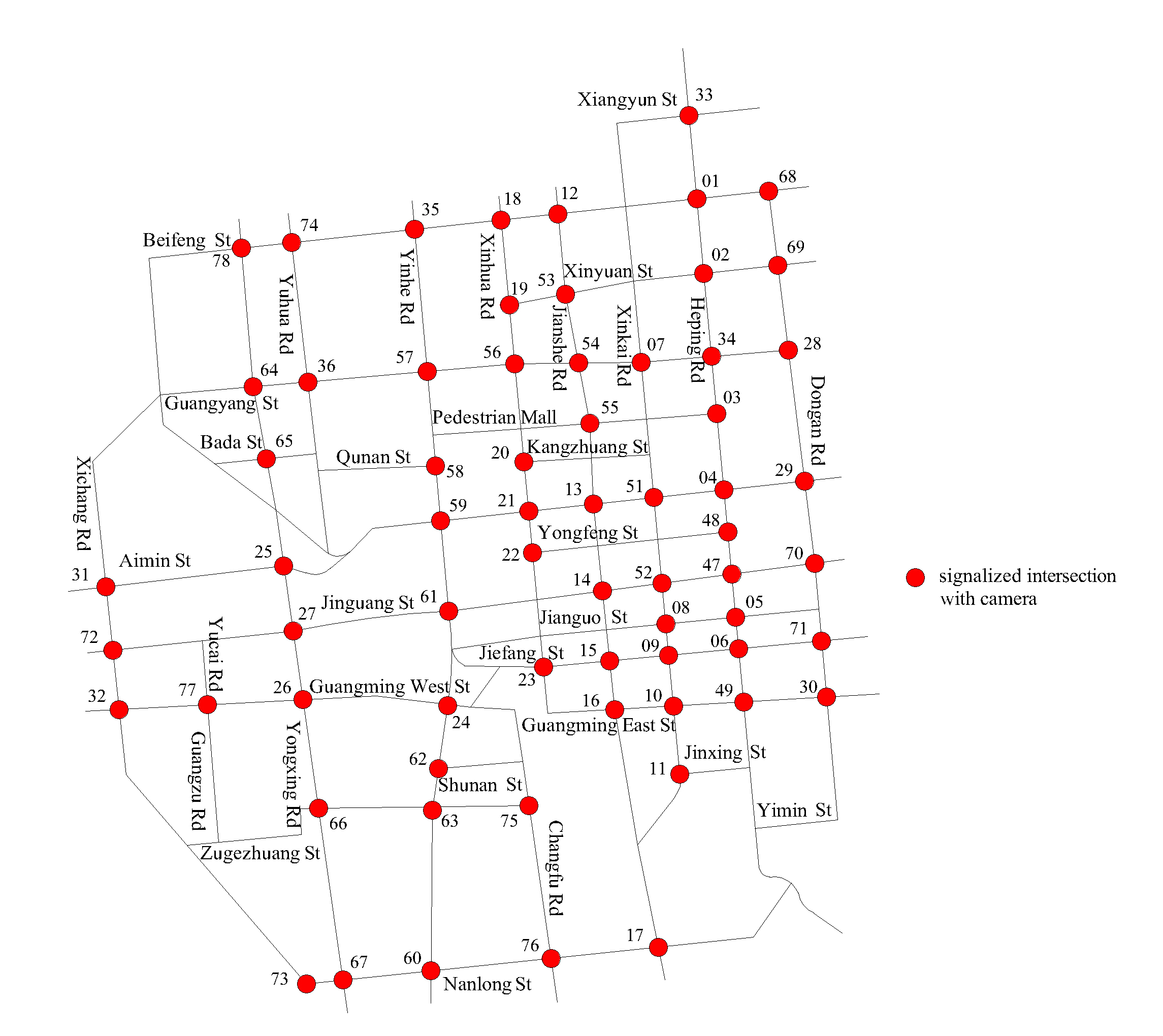}
  \caption{Data Location}\label{Fig:loc}
\end{figure}

(II) We simplify the road network to a square grid of $10\times 10$  considering the scale ($\sim 7km \times 7 km$) and the shape (Figure \ref{Fig:loc}) of the city. The original nodes are mapped into nodes in the grid with proper rescaling and rotating.

(III) We use a triangular flow-density model \cite{Ni_2015} for the speed-density dynamics. Specifically, given the free flow density limit $k_{m}$, congestion density limit $k_{c}$, and free flow speed $v_{m}$, the flow-density relationship between flow $q$ and density $k$ can be described as follows:
\begin{equation}\nonumber
  q =
  \begin{cases}
    v_{m}k & k<k_m\\
    \frac{k_c-k}{k_c-k_m}v_{m}k_m & k_m<k<k_c
  \end{cases}
\end{equation}

(IV) We use a linear function to combine the traffic density in the simulation with the background traffic density. For time $t$ and link $l$, the estimated travel density $\hat{k}_{tl}$ is calculated as
\begin{equation}\nonumber
  \hat{k}_{tl} = (1-\varphi)k_{tl} + \frac{1}{k_0}k'_{tl}
\end{equation}

\noindent where $k_{tl}$ is the density derived from the background travel time, $k'_{tl}$ is the density of vehicle/trips in the simulator, and $\varphi$ is the market exposure rate of the operator's service, which quantifies the share of the population that can be represented by the simulation process. Conversion factor $k_0$ is a parameter to ensure that $k_{tl}'$ and $k_{tl}$ have the same scale, as in simulation we generally model much fewer vehicles to reduce computation requirement. We calibrate this factor to make the background travel time distribution consistent across different market exposure rate $\varphi$. The value of $k_0$ for specific setting is reported along with the result in next section.

\begin{figure}[!t]
  \centering
  \includegraphics[width = 0.5\textwidth]{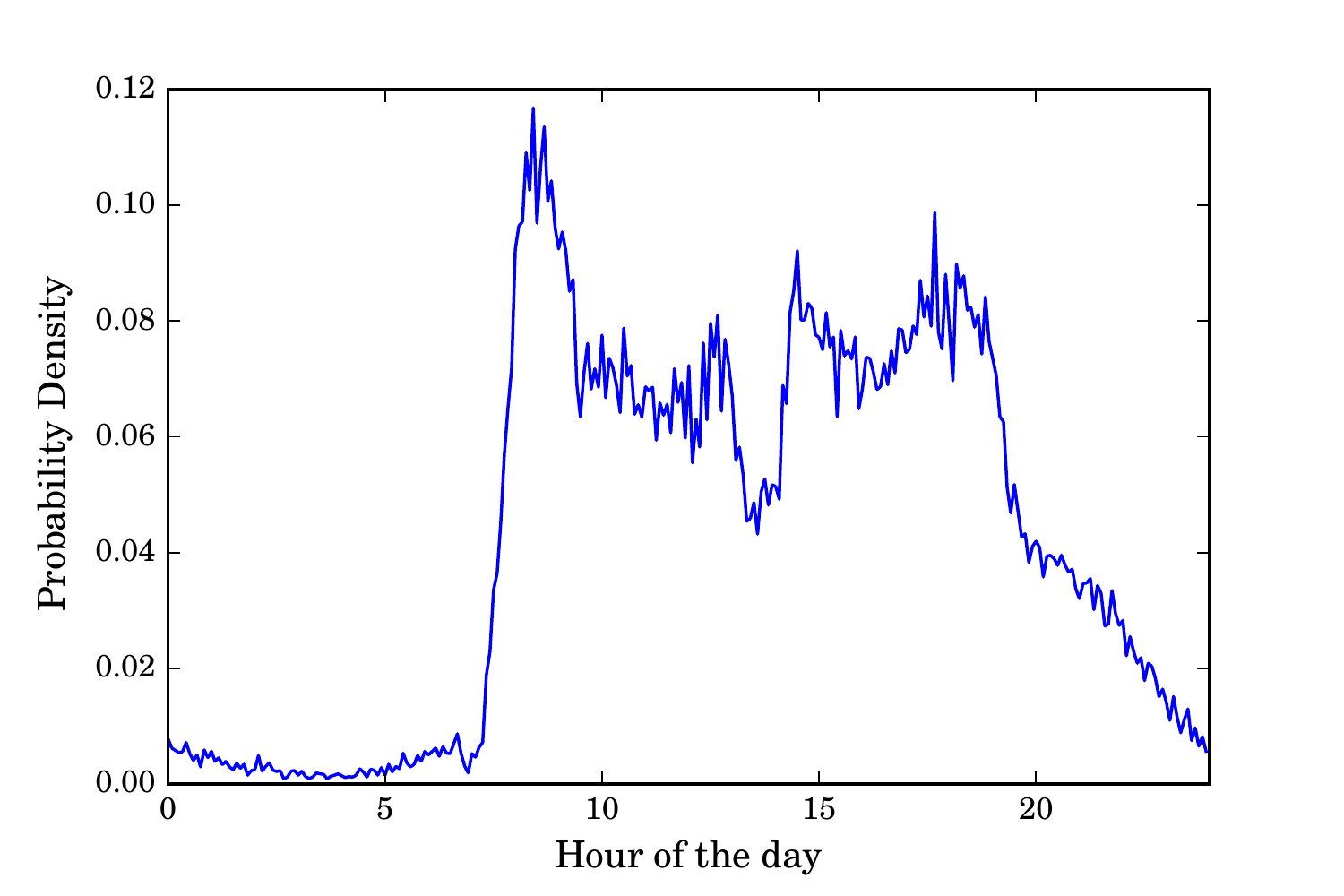}
  \caption{Distribution of Trips across Time}\label{Fig:demand_ts}
\end{figure}

As introduced above, the simulation system takes background traffic density and the distribution of travel demand as inputs. We obtain these input from our dataset, which consists of real trajectory data of private trips in the Langfang central region on a normal day. There are about 250,000 trip records in total, and each record is a tuple consisting of series of location identified code and series of time stamps in seconds, i.e., $(\{l_1,\cdots,l_n\},\{t_1,\cdots,t_n\})$. Figure \ref{Fig:demand_ts} shows the distribution of trips across time.

To obtain the demand distributions, we extract the origin and destination information from the trajectory records and then aggregate them for each time step. Such information is then normalized and smoothened to form demand distributions for each OD pair between cells in the grid at each time step. We reserve a degree of freedom by defining a rescaling factor $M$ to control the total number of trips in a day, which enable us to balance between reducing computation burden and controlling variances in the results.

To obtain the background travel density, we simply extract and aggregate the link-level information at each time step from trajectories. To deal with possible missing elements, we apply interpolation and extrapolation from nearby observations. Finally, a uniform filter is applied to smoothen the density distribution.

Next, we describe the parameters we selected for the experiment.

\subsubsection{MoD Service Operation}

For operating cost, we only consider cost of fuel (in following linear form)
\begin{equation}\nonumber
  c = c_d\cdot d
\end{equation}

\noindent where $c_d$ is the unit distance cost and set to be $\$0.07/km$ based on the current fuel costs ($\$1/L$) \cite{Data_Gas} and standard fuel efficiency ($7L/100km$) \cite{Data_Fuel_Eff} in China.

Both the fare for single services $f_s$ and that for shared services $f_{sh}$ are related to distance $d$ and time $T$
\begin{equation}\nonumber
  \begin{aligned}
    f_s & = r_b + r_d\cdot d + r_t\cdot T\\
    f_{sh} & = (r_b + r_d\cdot d + r_t\cdot T)r_{sh}
  \end{aligned}
\end{equation}

\noindent where $r_b,r_d,r_t$ are the base fare, the unit distance fare, and the unit time fare  for each trip and are set to be $\$1$, $\$0.25/km$, and $\$0.01/min$ respectively based on the taxi fare structure in Langfang \cite{Data_Taxi}. $r_{sh}(<1)$ is the uniform fare ratio between shared services and single services and varies in different scenario settings.

For operation, we have following rules. A maximum of three parties are allowed to share one vehicle. For each request, the operator searches for empty vehicle within a maximum pickup distance of $5 \ km$ and non-empty vehicles with additional capacity within a maximum pickup distance of $2 \ km$. We also restrict the maximum detour distance and time of shared trips, which are defined as the difference between distance and travel time of shared service and those of single service, to $2 \ km$ and $5 \ mins$ respectively.

\subsubsection{Specification for Traveler Behavior}

We follow \cite{Atasoy_2015} to specify a choice model considering the service type, travel mode, travel time, and travel cost. We define the utility value in the monetary unit for a realistic selection of parameters. The detailed specification is shown below, where $U_s,U_{sh},U_o$ and $U_n$ are the utilities of taking the single service, the shared service, the original travel mode, or cancellation.
\begin{equation}\nonumber
  \begin{aligned}
    & U_s = \mu(ASC_s - VOT\cdot T_s - f_s - E(\delta_s)) \\
    & U_{sh} = \mu(ASC_{sh} - b_{v,sh} VOT\cdot T_{sh} - f_{sh} - E(\delta_{sh})) \\
    & U_o = \mu(ASC_o - VOT\cdot T_o - b_{f,o}f_o) \\
    & U_n = 0
  \end{aligned}
\end{equation}

In the specification above, factor $VOT$ characterizes the average value of time of this population and is set to be $\$0.03/min$, $70\%$ of the average income ($\$2.5/hour$) in Langfang \cite{Data_Salary,Wu_2014,Data_lowest_income}.
Factor $\mu$ is the scale factor and set to be $0.5$ as in the work of Atasoy \textit{et al.} \cite{Atasoy_2015}.
The constant $ASC$s characterize the utility not captured by travel time and cost, and are set as $ASC_s = 4.5, ASC_{sh} = 4, ASC_o = 5$.

$E$ reflects the asymmetry in perception of utility between upward and downward fare adjustment, and it is defined to be a piecewise linear function
\begin{equation}\nonumber
  E(\delta) =
    \begin{cases}
      e_1\delta & \delta < 0\\
      e_2\delta & \delta \geqslant 0
    \end{cases}
\end{equation}

\noindent with $e_2 = 2$, $e_1 = 1$. This function reflects the fact that people are more sensitive towards price surges than price reductions \cite{Kahneman_1979}. Specifically, in our simulation we assume travelers regard price reductions as normal fare changes, whereas price surges are considered as unexpected and generate twice disutility.

Factor $b_{v,sh}$ characterizes the relative time reliability of shared services compared with single services and is define to be $\eta_{sh} / \eta_s$, where $\eta_{sh}$ is the ratio between the actual travel time and the estimated travel time for shared services, and $\eta_s$ is the ratio for single services. This ratio is set to be $1.2$ based on our simulation results in the context of Langfang.

Factor $b_{f,o}$ characterizes the average cost of the original travel mode compared with driving and varies by individuals and their car ownership, transit accessibility, etc. We control this factor to specify different settings of attractiveness of other alternatives for the population. A lower $b_{f,o}$ corresponds to higher car ownership and lower cost of staying in the original travel mode across the population, and launching the new service has less impact; whereas a higher $b_{f,o}$ represents a setting with lower car ownership, the original travel mode is costly and the new service can have a greater impact. This factor is determined in specific scenario setting.

\subsubsection{Flow-density Model}

Based on Langfang's local traffic situation we choose $v_m = 60km/h$, $k_m = 1veh/36m$ and $k_c = 1veh/6m$. Since the scale of $k$ is irrelevant if we only need the information of travel time, in the implementation we use $k_m = 1$ and $k_c = 6$ to simplify the computation.

\subsection{Results}

We implement this simulation system in C++, and simulate each scenario on a separate Linux server in the cluster. Each server has 64GB memory and a 8-core 2.0GHz CPU. In each step of CMA-ES algorithm, 12 different $\theta$ are sampled and for each the function $J(\pi_{\theta})$ is evaluated and averaged over 10 runs. In our experiment, with the number of total request $M=100,000$, each step takes about $15 \ mins$, and we allow the the algorithm to run for 20 steps ($5 \ hour$).

Our algorithm is tested on several scenario settings, with different combinations of fleet size $N$, shared service fare ratio $r_{sh}$, original cost factor $b_{f,o}$, market exposure $\varphi$, and background congestion levels $\psi_C$.

To reduce the variance in aggregate results from simulation, we set the potential request size to $M = 100,000$. We then consider $r_{sh}\in\{0.4,0.6,0.8\}$ and $N\in\{500,750,1000,1250\}$. We select these $N$ based on the observation that for $N \geqslant 1250$, the increase in the number of passengers with greater $N$ becomes negligible.

We consider $b_{f,o}=2.5,5.0$ to model the trip generation ability of this MoD service (comparing with original mode). We set the density adjustment factor $k_0 = 0.95,0.8$ correspondingly to keep density combination process consistent.

$\varphi$ and $\psi_C$ are used to capture the potential impact of the private MoD operation on system congestion level and have the following range $\varphi\in\{10\%,20\%,40\%\},\psi_C\in\{\text{Low, Medium, High}\}$ respectively.
Each level of $\psi_C$ is derived by multiplying the background travel density of each link by a certain factor: the low, medium, and high level of $\psi_C$ correspond to $\text{factor} = 0.8,1.0,1.2$, respectively.

In comparison, we consider performance of basic strategies, the myopic pricing strategy (PM) and the proposed optimal pricing strategy (PO). Basic strategies include the one to provide all feasible options if possible (S+Sh), to only provide feasible single options (S), and to only provide feasible shared options (Sh).
We include several metrics to assess performance on different aspects:

(I) the profit, revenue, and cost of the operations to assess the financial performance of the private MoD operator;

(II) the number of passengers serviced, and the total travel time and distance of the whole fleet to assess the operation efficiency;

(III) the average waiting time (AWT) for pickup across serviced passengers, which is defined to be gap between the time of the trip confirmation and the time of actual pickup, to assess the service quality;

(IV) the global capacity, defined to be the ratio $\rho$ between the number of trips realized (including single service, shared service, and original travel modes), and the number of total requests:
\begin{equation}\nonumber
    \rho = \frac{M_s + M_{sh} + M_o}{M_s + M_{sh} + M_o + M_n}
\end{equation}

\noindent where $M_s$, $M_{sh}$, $M_o$, and $M_n$ are the number of single trips, shared trips, original trips, and trip cancellations from the target population, respectively;

(V) the congestion level, defined to be the average system unit link travel time
\begin{equation}\nonumber
  T_{ave} = \frac{\sum_t\sum_{e\in links} k_{e,t}T_{e,t}}{\sum_t\sum_{e\in links} k_{e,t}}
\end{equation}

\noindent where $T_{e,t}$ and $k_{e,t}$ are the unit link travel time and density for link $e$ at time $t$, respectively;

(VI) the efficiency of sharing, measured by total reduced trip distance $\Delta_d$, or the difference between the sum of (direct) trip distance of realized trips and the total traveled distance by vehicles
\begin{equation}\nonumber
  \Delta_d = \sum_{trip\in trips}d_{trip} - \sum_{veh\in fleet}d_{veh}.
\end{equation}

The comparisons on operating metrics (I-III) are presented as the percentage difference for selected strategies, compared with the optimal pricing strategy (PO), averaged over 108 different settings of $(N,r_{sh},\rho,\psi_C)$, under different $b_{f,o}$ in Table \ref{tab:statistics_dyna}. The optimal pricing strategy consistently generates a much higher profit than the other strategies, as expected; but this comes at the cost of serving fewer trips. Interestingly, for PO, the ratio between the total travel time and the total travel distance of the fleet is smaller than that of others, which means that PO attempts to avoid trips that pass through a congested area. In terms of service quality, for PO, the average waiting time for pickup is much shorter than that of others.

\begin{table}[!t]
\centering
\caption{Average percentage difference in operating statistics, compared with PO}\label{tab:statistics_dyna}
\begin{tabular}{cccccc}
\hline
\multicolumn{5}{c}{$b_{f,o} = 2.5$} \\ \hline
Strategy & S+Sh & S & Sh & PM\\ \hline
Profit & -11.18\% & -39.51\% & -26.10\% & -11.45\% \\
Revenue & -7.99\% & -39.32\% & -22.61\% & -8.22\% \\
Cost & 7.80\% & -37.81\% & -5.61\% & 8.42\% \\
Serviced Passengers & 22.09\% & -44.14\% & 15.40\% & 16.78\% \\
Travel Distance & 7.80\% & -37.81\% & -5.61\% & 8.42\% \\
Travel Time & 12.40\% & -30.54\% & -2.54\% & 12.16\% \\
AWT, Pickup & 13.65\% & -8.39\% & 10.70\% & 13.31\% \\\hline
\multicolumn{5}{c}{$b_{f,o} = 5.0$} \\ \hline
Strategy & S+Sh & S & Sh & PM\\ \hline
Profit & -16.80\% & -36.01\% & -25.52\% & -14.90\% \\
Revenue & -13.01\% & -33.98\% & -21.65\% & -11.75\% \\
Cost & 7.58\% & -21.96\% & -0.94\% & 5.86\% \\
Serviced Passengers & 25.51\% & -35.68\% & 24.03\% & 17.42\% \\
Travel Distance & 7.58\% & -21.96\% & -0.94\% & 5.86\% \\
Travel Time & 14.05\% & -13.31\% & 4.01\% & 12.07\% \\
AWT, Pickup & 17.71\% & 6.03\% & 14.01\% & 17.94\% \\\hline
\end{tabular}
\end{table}

\begin{table}[!t]
\centering
\caption{Background system metrics}\label{tab:background}
\begin{tabular}{ccc}
\hline
\multicolumn{3}{c}{$b_{f,o} = 2.5$} \\ \hline
$\psi_C$ & $T_{ave}$($mins/km$) & $\rho$ \\ \hline
Low & 1.666 &0.844\\
Medium &2.272 &0.832\\
High &3.063 &0.815\\\hline
\multicolumn{3}{c}{$b_{f,o} = 5.0$} \\ \hline
$\psi_C$ & $T_{ave}$($mins/km$) & $\rho$ \\ \hline
Low & 1.662 &0.745\\
Medium &2.277 &0.725\\
High &3.088 &0.699\\\hline
\end{tabular}
\end{table}

\begin{figure}[!t]
  \centering
  \subfigure[$b_{f,o} = 2.5$]{
    \includegraphics[width = 0.475\textwidth]{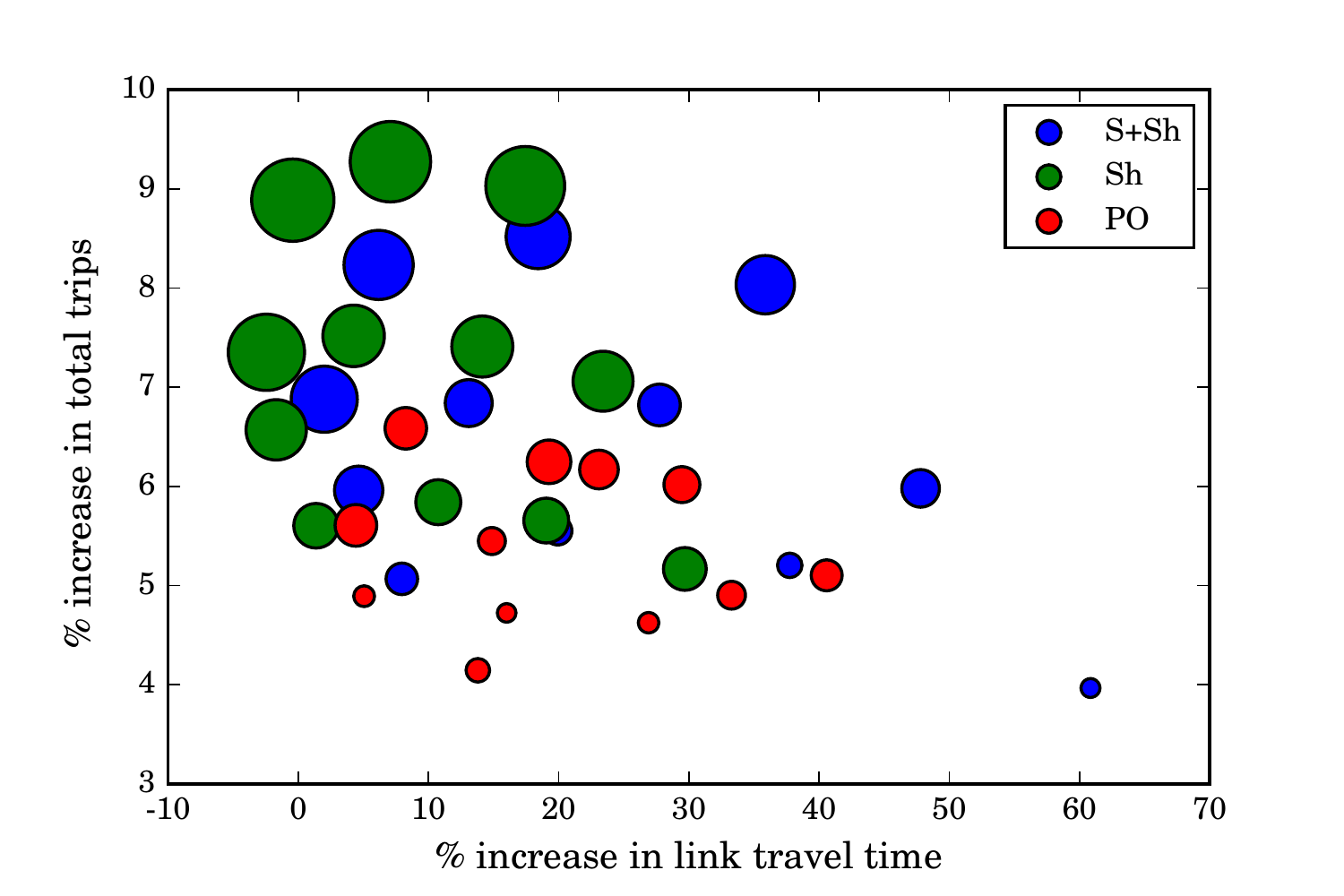}
  }
  \subfigure[$b_{f,o} = 5.0$]{
    \includegraphics[width = 0.475\textwidth]{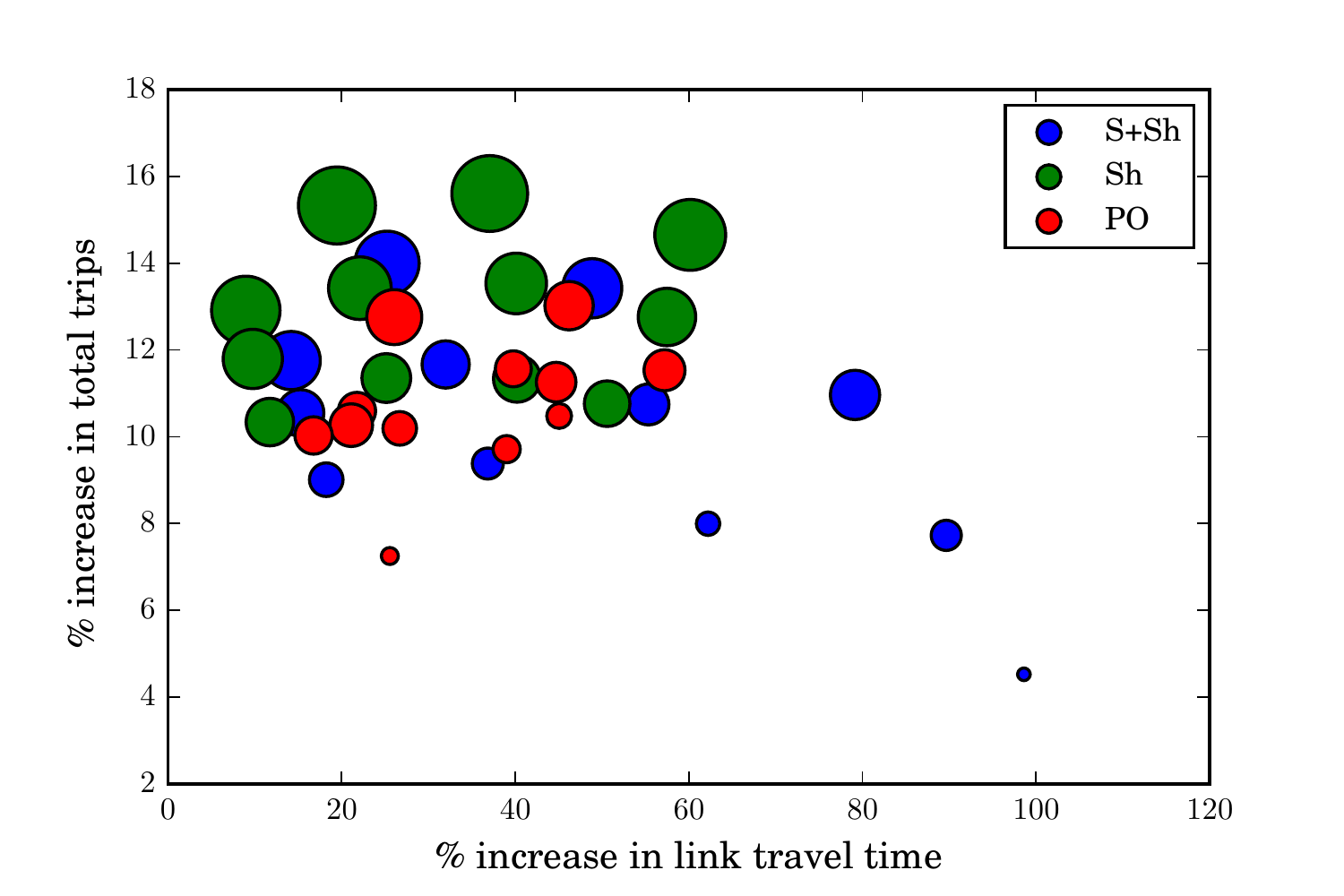}
  }
  \caption{Congestion - Capacity Tradeoff among Strategies}\label{fig:cc_dyna}
\end{figure}

The comparisons on system level metrics (IV-VI) are presented in Figure \ref{fig:cc_dyna}.
In this figure, circles correspond to different strategies and settings of $(N,r_{sh})$ under high market share and high background congestion $\rho=40\%$ and $\psi_C=\text{High}$. Each circle has its radius proportional to $\Delta_d$, and its center pointed by the percentage differences in $T_{ave}$ and $\rho$ compared with the background value in Table \ref{tab:background} (which referred to system without the impact of the design MoD service. We only consider S+Sh and Sh as benchmarks since they are more likely to have better system performances, as shown in the Table \ref{tab:statistics_dyna}.

Figure \ref{fig:cc_dyna} shows that, compared with those of other two, points of PO are generally located in the bottom right side and have smaller radii, indicating the system-level inefficiency of PO. But these points are more closely clustered, implying that dynamic pricing introduces robustness towards changes in settings, such as those under severe congestion.

\section{Conclusion}

We consider a profit maximization problem in shared mobility on-demand service with dynamic pricing decisions, and design an efficient algorithm to solve for near-optimal strategies. Based on a simulated experiment in Langfang, China, we show that the proposed optimal pricing strategy (PO) has superior performance in profit compared with other strategies. However, this strategy can lead to inefficiency of the system, mainly because it has stricter criteria for providing shared capacity. The conflict between the profit maximization objective and the system capacity and efficiency objective highlights the need to understand the perspective of the private operator and potentially introduce policy interventions to better align the two objectives.

There are many potential extensions of this work. First, in our problem formulation, we assume that the operator has full distributional knowledge of the demand. This assumption is too strong in practice and a practical extension is to consider simultaneous online estimation of behavior model and behavior-model-based profit maximization \cite{Rusmevichientong_2010,Saure_2013,Caro_2007}.

Another extension is to integrate the pricing problem with the dispatching problem. The optimal strategies of these two problems should be closely related; in fact, in this work, our optimal pricing strategy show the ability to allocate supply and incentivize demand wisely. Thus we expect that the integrated problem can be formulated and solved under similar non-myopic and reinforced framework.

% An example of a double column floating figure using two subfigures.
% (The subfig.sty package must be loaded for this to work.)
% The subfigure \label commands are set within each subfloat command,
% and the \label for the overall figure must come after \caption.
% \hfil is used as a separator to get equal spacing.
% Watch out that the combined width of all the subfigures on a
% line do not exceed the text width or a line break will occur.
%
%\begin{figure*}[!t]
%\centering
%\subfloat[Case I]{\includegraphics[width=2.5in]{box}%
%\label{fig_first_case}}
%\hfil
%\subfloat[Case II]{\includegraphics[width=2.5in]{box}%
%\label{fig_second_case}}
%\caption{Simulation results for the network.}
%\label{fig_sim}
%\end{figure*}

\appendices

\section{Proof of Proposition 1}

1) From \eqref{p_p} to \eqref{e_p}:

Assume the optimal value of \eqref{p_p} is $z^*$. By the optimality of $z^*$, we know that for any $\delta\in R^n$
\begin{equation}\nonumber
  \begin{aligned}
  & z^* \geqslant \frac{\sum_{j=1}^n e^{U_j(\delta_j)}(f_j-c_j+\delta_j)}{\sum_{j=1}^n e^{U_j(\delta_j)} + e^{U_0}} \\
  \Rightarrow & e^{U_0}z^* \geqslant \sum_{j=1}^n e^{U_j(\delta_j)}(f_j - c_j - z^* + \delta_j)
  \end{aligned}
\end{equation}

But there is also some $\delta^*$ that reach optimality for \eqref{p_p}, so we have
\begin{equation}\nonumber
  e^{U_0}z^* = \max_{\delta\in R^n}\sum_{j=1}^n e^{U_j(\delta_j)}(f_j - c_j - z^* + \delta_j)
\end{equation}

\noindent which is exactly \eqref{e_p}. Thus, $z^*$ is a solution of \eqref{e_p}.

2) From \eqref{e_p} to \eqref{p_p}:

If $z^*$ solves \eqref{e_p}, we know that for every $\delta \in R^n$,
\begin{equation}\nonumber
  \begin{aligned}
  & e^{U_0}z^* \geqslant \sum_{j=1}^n e^{U_j(\delta_j)}(f_j - c_j - z^* + \delta_j)\\
  \Rightarrow & z^* \geqslant \frac{\sum_{j=1}^n e^{U_j(\delta_j)}(f_j-c_j+\delta_j)}{\sum_{j=1}^n e^{U_j(\delta_j)} + e^{U_0}},
  \end{aligned}
\end{equation}

\noindent thus $z^*$ is an upper bound of \eqref{p_p}. But there is also some $\delta^*$ that reach optimality for the right-hand side of \eqref{e_p}, so we must have $z^*$ to be the strict upper bound of \eqref{p_p}.

\section{Proof of Proposition 2}

By proposition 1, we only need to prove the process converges into the solution of \eqref{e_p}. We note that the left-hand side of \eqref{e_p} is a monotonically increasing linear function of $z$ on $R$, whereas the right-hand side is a monotonically non-increasing convex function of $z$ on $R$. The later claim is valid as the right-hand side is the element-wise maximum of a collection of convex function
\begin{equation}\nonumber
  h_{\delta}(z) = \sum_{j=1}^n e^{U_j(\delta_j)}(f_j - c_j - z + \delta_j)
\end{equation}

Thus, $h(z) = \max_{\delta\in R^n} h_{\delta}(z) - e^{U_0}z$ is convex and monotonically decreasing on $R$. Moreover, this function is positive at $z=0$, as $\max_{\delta\in R^n} h_{\delta}(z)$ is positive for all $z$ by
\begin{equation}\nonumber
  \max_{\delta\in R^n}h_{\delta}(z) \geqslant \sum_{j=1}^n e^{U_j(z+c_j-f_j+\varepsilon)}\varepsilon > 0
\end{equation}

\noindent with some $\varepsilon > 0$. Therefore, the zero of $h(z)$ at $R_+$ is unique and can be solved by Newton method based on the monotonicity and convexity of $h(z)$
\begin{equation}\nonumber
  \begin{aligned}
  z^{k+1} & = z^k - h(z^k) / h'(z^k) \\
  & = z^k - \frac{\sum_{j=1}^n e^{U_j(\delta_j^k)}(f_j - c_j - z^k + \delta_j^k) - e^{U_0}z^k}{- \sum_{j=1}^n e^{U_j(\delta_j^k)} - e^{U_0}} \\
  & = \frac{\sum_{j=1}^n e^{U_j(\delta_j^k)}(f_j - c_j + \delta_j^k)}{\sum_{j=1}^n e^{U_j(\delta_j^k)} + e^{U_0}}
  \end{aligned}
\end{equation}

\noindent where $\delta^k = \arg\max_{\delta\in R^n}h_{\delta}(z^k)$. This is exactly \eqref{re_p:2} if we expand $e^{U_j(\delta_j^k)}$ to $V_{j0}e^{E_{y_j}(\delta_j^k)}$. We therefore remain to show that $\delta^k$ is of form \eqref{re_p:1} under our assumption on the form of $E_{y_j}$ as
\begin{equation}\label{lem:1}
  \delta_j^k = \arg\max_{\delta_j\in R} e^{E_{y_j}(\delta_j)}(f_j - c_j - z^k + \delta_j)
\end{equation}

When $E_{y}(x) = - e_{y1}\min\{0,x\} - e_{y2}\max\{0,x\}$ with $e_{y2} > e_{y1} > 0$ for all service type $y$, the derivatives of the function in \eqref{lem:1} can be calculated separately on $\delta_j > 0$ and $\delta_j < 0$, and is equal to
\begin{equation}\nonumber
    \begin{cases}
      e_{y_j1}e^{-e_{y_j1}\delta_j}(z^k-f_j+c_j+1/e_{y_j1}-\delta_j) & \delta_j < 0\\
      e_{y_j2}e^{-e_{y_j2}\delta_j}(z^k-f_j+c_j+1/e_{y_j2}-\delta_j) & \delta_j > 0\\
      undefined & \delta_j = 0
    \end{cases}
\end{equation}

Since $1/e_{y_j1} > 1/e_{y_j2}$, there will be at most one stationary point $\delta_j$. As the target function in \eqref{lem:1} is second order differentiable in $R/\{0\}$, the optimal $\delta_j^k$ is the stationary point if there is one, or is 0 if there is no stationary point, and the derivative is positive when $\delta_j < 0$ and negative when $\delta_j > 0$. This is exactly the cases given by \eqref{re_p:1}.

% use section* for acknowledgment
\section*{Acknowledgment}

% Can use something like this to put references on a page
% by themselves when using endfloat and the captionsoff option.
\ifCLASSOPTIONcaptionsoff
  \newpage
\fi

% trigger a \newpage just before the given reference
% number - used to balance the columns on the last page
% adjust value as needed - may need to be readjusted if
% the document is modified later
%\IEEEtriggeratref{8}
% The "triggered" command can be changed if desired:
%\IEEEtriggercmd{\enlargethispage{-5in}}

% references section

% can use a bibliography generated by BibTeX as a .bbl file
% BibTeX documentation can be easily obtained at:
% http://mirror.ctan.org/biblio/bibtex/contrib/doc/
% The IEEEtran BibTeX style support page is at:
% http://www.michaelshell.org/tex/ieeetran/bibtex/
%\bibliographystyle{IEEEtran}
% argument is your BibTeX string definitions and bibliography database(s)
%\bibliography{IEEEabrv,../bib/paper}
%
% <OR> manually copy in the resultant .bbl file
% set second argument of \begin to the number of references
% (used to reserve space for the reference number labels box)
\bibliographystyle{IEEEtran}
\bibliography{IEEEabrv,dyna_pricing}

% biography section
%
% If you have an EPS/PDF photo (graphicx package needed) extra braces are
% needed around the contents of the optional argument to biography to prevent
% the LaTeX parser from getting confused when it sees the complicated
% \includegraphics command within an optional argument. (You could create
% your own custom macro containing the \includegraphics command to make things
% simpler here.)

% or if you just want to reserve a space for a photo:

\begin{IEEEbiographynophoto}{Han Qiu}
Han Qiu received the B.E. and B.S. degrees from Tsinghua University, Beijing, China, in 2015, and the M.S. degree from Massachusetts Institute of Technology, Massachusetts, USA, in 2017. His main research interests include human behavior and reinforcement learning.
\end{IEEEbiographynophoto}

\begin{IEEEbiographynophoto}{Ruimin Li}
Ruimin Li is an associate professor of Institute of transportation engineering in the Department of Civil Engineering at Tsinghua University, China. He received his B.S. and Ph. D from the Tsinghua University in 2000 and 2005. He was born in 1979 in Shandong province, China. He joined the Institute of Transportation Engineering, Tsinghua University in 2005. His current research interests include intelligent transportation system, traffic control, traffic management, traffic simulation.
\end{IEEEbiographynophoto}

\begin{IEEEbiographynophoto}{Jinhua Zhao}
Jinhua Zhao is the Edward and Joyce Linde Associate Professor of City and Transportation Planning at the Massachusetts Institute of Technology (MIT). Prof. Zhao brings behavioral science and transportation technology together to shape travel behavior, design mobility system and reform urban policies. He develops methods to sense, predict, nudge and regulate travel behavior, and designs multimodal mobility system that integrates autonomous vehicles, shared mobility and public transport. Prof. Zhao directs the Urban Mobility Lab (mobility.mit.edu) at MIT.
\end{IEEEbiographynophoto}

\end{document}